\documentclass[12pt]{article}
\usepackage{amsmath,amsthm,amssymb}

\newtheorem{theorem}{Theorem}

\newtheorem{lemma}{Lemma}

\newtheorem{corollary}{Corollary}

\newcommand{\ga}{\alpha}
 
\newcommand{\gc}{\gamma}
\newcommand{\gd}{\delta}
\newcommand{\gl}{\lambda}

\newcommand{\gra}{\nabla}
\newcommand{\de}{\partial}
\newcommand{\bpf}{\begin{proof}}
\newcommand{\epf}{\end{proof}}
\newcommand{\beq}{\begin{equation}}
\newcommand{\eeq}{\end{equation}}

\begin{document}
\title{Asymptotic Behaviors of Non-variational Elliptic Systems} 
 \author{Szu-yu Sophie Chen
  \footnote{The author was supported by the Miller Institute for Basic Research in Science.}}
 \date{}
\maketitle

 \begin{abstract}
We  use a  method, inspired by Pohozeav's work, 
to study asymptotic behaviors of non-variational elliptic systems in dimension $n \geq 3.$  
As an application, we prove removal of an apparent singularity in a ball and uniqueness 
 of the entire solution. All results apply to changing sign solutions.
\end{abstract}

 A classical work by   Gidas and Spruck \cite{GS81} asserts that any nonnegative solution to 
  $\Delta u + |u|^{\ga-2} u = 0$  in $\mathbb R^n$ with $2 < \ga < \frac{2n}{n-2}$ (subcritical case) is trivial.
  For $\ga = \frac{2n}{n-2},$ Caffarelli, Gidas and Spruck \cite{cgs89} proved that any nonnegative solution in $\mathbb R^n$ is of the form $u = (a+ b |x|^2)^{- \frac{n-2}{2}},$
  where $a, b$ are constants. Such problem for elliptic systems are also studied,  
  for example,  in the studies of Lane-Emden type systems; see  Zou \cite{Zou} and Polacik, Quittner and Souplet \cite{PQS}
   and the references therein. 
 
  On the other hand, the behaviors of changing sign solutions are more delicate. For example, there exists a sequence 
  of changing sign solutions to  $\Delta u + |u|^{\ga-2} u = 0$  in $\mathbb R^n$ with $2 < \ga < \frac{2n}{n-2};$ see \cite{KP97}.
  In this paper,  we study under what circumstances a solution  to an  elliptic system in an exterior domain is asymptotic to  $|x|^{-(n-2)}$ at the
 infinity.  Such  decay  is optimal in the sense that  the infinity is a regular point in the inverted  coordinates.
   It is known  \cite{KP97} that  there exist  solutions  to $\Delta u + u^{\ga-1} = 0$ in $\mathbb R^n$ which decay slower than $|x|^{-(n-2)}.$ Thus, a suitable integrability condition is necessary to exclude such case.

  While the study of  changing sign solutions to elliptic systems is interesting by itself, the problem is  well-motivated by differential geometry.
  For example, the decay of curvature tensors was studied for Yang-Mills fields \cite{Uh82}, Einstein metrics \cite{BKN89}
   and other generalizations   \cite{TV05}, \cite{Chen08}, just to name a few.  A typical system for the problem is of the  form
    $$        \Delta (Rm)_{ijkl} =  Q_{ijkl}(Rm, Rm),       
   $$ where $Rm$ is the Riemannian curvature tensor and $Q$ is a quadratic in $Rm.$  A natural geometric assumption is that 
    $|Rm|$ is in $L^{\frac{n}{2}}.$ Therefore, $|Rm|$ vanishes at infinity and the problem is to find out the decay rate.   
    The study of geometrical systems is more subtle as $(Rm)_{ijkl}$ satisfies an extra relation,
    the Bianchi identity, and  the underlying spaces are not Euclidean.  
   
 The technique we use in this paper is based on the method developed in \cite{Chen08}  on asymptotically flat manifolds, where
  a special geometric setting is considered.  In this paper, we study general non-variational elliptic systems of the reaction-diffusion type. 
 Our result applies to changing sign solutions and includes the supercritical case (i.e.,  $\Delta u + C u^{\ga-1}  = 0$  with $ \ga > \frac{2n}{n-2},$ where $C$  is a constant).

Let $V = (V_1, \cdots, V_m)$ and $f^i: \mathbb{R}^m \rightarrow  \mathbb{R}.$  
Consider the system of equations 
\beq \label{e:system}
 \sum_{j=1}^m  A_{ij} \Delta V_j =   f^i(V),
  \eeq where $A$ is a constant  positive  definite symmetric matrix and  $i = 1, \cdots, m.$
  The system (\ref{e:system})  describes the steady states of the reaction-diffusion systems.
   The matrix $ A_{ij} $ represents the diffusion rate  and $f^i(V)$ is the reaction term. 
    Note that   $f^i(V)$ or $V_i$  may have no  sign. 
   We assume the following structure conditions
  
    (A1)  $| f^i(V)| \leq C |V|^q,$ 
          
  (A2)  $| \gra f^i(V)| \leq C |V|^{q-1}.$  
  
  Let $K$ be a compact subset in $\mathbb{R}^n.$
 \begin{theorem} \label{t:asymp} Let $q > \frac{n+2}{n}$ and $p = \frac{n}{2} (q-1).$
  Suppose that $f^i$ satisfies (A1) and (A2).
  Let  $V \in L^p(\mathbb{R}^n \setminus  K)$ be a solution to (\ref{e:system})  in  $\mathbb{R}^n \setminus K.$
   Then $|V| = O(|x|^{- (n-2)})$  and $|\gra V| = O(|x|^{- (n-1)})$ at the infinity.
  \end{theorem}
 
  An immediate consequence is a result on singularity removal  for  affine invariant equations.
 For scalar equations, the problem was studied in \cite{GS81}, \cite{BL81}, \cite{cgs89}. 
 
 Let $B_1$ be the unit ball centered at the origin. 
 \begin{corollary}  \label{t:removal}
Suppose $f^i$ are homogeneous functions of degree $\frac{n+2}{n-2}.$
  Let $V \in L^{\frac{2n}{n-2}}(B_1)$ be a solution to (\ref{e:system}) in $B_1\setminus  \{0\}.$
   Then $V$ can be extended to a smooth solution to (\ref{e:system}) in  $B_1.$
\end{corollary}

  By performing a linear transformation $W_i = \sum_j A_{ij} V_j,$ the system (\ref{e:system}) can be reduced to 
  an equation of  the diagonal form $\Delta W = \tilde f(W).$  
  The assumptions (A1)-(A2) and other conditions on $V$ or $f^i$  equivalently hold for  $W$ and $\tilde f.$
  Therefore, for Theorem~\ref{t:asymp} and Corollary~\ref{t:removal}, we may assume without loss of generality 
  the equation is of  the diagonal form.
 
   We turn to study the uniqueness of entire solutions for variational systems.
   Let $P(V)$ be a homogeneous function of degree $q+1.$ Suppose that $f^i = \frac{\de P}{\de V^i}$ in 
    (\ref{e:system}). Note that (\ref{e:system}) can not be reduced to a variational system of the form $\Delta W = \tilde f(W)$
    now; it can still be reduced to one with diagonal matrix $A_{ij}.$
      For scalar equations,   there is a large literature on the uniqueness problem; e.g. Gidas and Spruck \cite{GS81}, Bidaut-Veron \cite{Bidaut-Veron} and Serrin and Zou \cite{SZ02}; see also \cite{PS07} and the references therein.
  For systems, when $P(V) \leq 0$ and $q > \frac{n+2}{n-2}$ (supercritical case),  the problem was studied by Pucci and Serrin \cite{PS86}  
 under some asymptotic assumption of $V$.  Their result also holds for non-homogeneous function $P$ (and more general $P(x, V, \gra V)$) satisfying some inequality.

\begin{theorem}  \label{t:unique}  Let $q > \frac{n+2}{n}$ $q \neq \frac{n+2}{n-2}$ and  $p = \frac{n}{2} (q-1).$
Suppose $P(V)$ is a homogeneous function of degree $q+1.$
   Let  $V \in L^p(\mathbb{R}^n)$ be a solution to (\ref{e:system}) in $\mathbb{R}^n$ with $f^i = \frac{\de P}{\de V^i}.$ 
     Then $V \equiv 0.$
\end{theorem}

We give the outline of proofs. To fix notations, 
  we denote by $dx$ the volume element in $\mathbb{R}^n$ and by $dS$ the area element of a hypersurface 
 in $\mathbb{R}^n.$
  Let $B_r (x)$ and $S_r(x)$ be the ball of radius $r$ and sphere of radius $r$ centered at $x,$ respectively. 
 When $x$ is at the origin, we simply denote by $B_r$ and $S_r.$ 

 The idea of  proof of Theorem~\ref{t:asymp} is to  compare the size of $\int_{\mathbb{R}^n \setminus B_r} |\gra V|^2 dx$ (as a function of $r$) to its derivative
  $-\int_{S_r} |\gra V|^2 dS.$ 
Then by ordinary differential inequality lemma,  we get the optimal decay of $|\gra V|$ and as a consequence the decay of $|V|.$
 In order to relate above two integrands, we use some version of  Pohozaev's identity for non-variational systems. 
   Pohozaev's  ingenious idea \cite{Poh65} is to use a conformal killing field to prove uniqueness in a star-shaped domain.
  This idea was generalized  nicely   by Pucci and Serrin \cite{PS86} to general variational systems. 
    Our use of the  identity is different 
  from the original one.  We apply  the identity to an unbounded domain (the complement of a large ball)  and use only the size of $|f^i|.$ 
 Therefore, our method can be applied to non-variational systems.
  
    The proof of Theorem~\ref{t:unique} is a combination of Theorem~\ref{t:asymp} and Pohozaev's original idea. 
    Since the solution decays fast enough at infinity, no terms from infinity contribute to the main integrand. 
    We use the identity differently such that  we obtain the uniqueness also in the subcritical case, in contrast to the problem
     in star-shaped regions where one has to restrict to the supercritical case.

  Finally, we show that the assumptions in above Theorems are sharp.
 
 Example 1: Consider the equation $\Delta u  + u^{\frac{n+2}{n-2}}= 0$ in $ \mathbb{R}^n.$  
      By  \cite{cgs89},  nonnegative solutions are of the form  $u= (a+ b|x|^2)^{- \frac{n-2}{2}}.$ Therefore, $u$ decays as $|x|^{-(n-2)}$
      at the infinity.   This example shows that in Theorem~\ref{t:unique}, 
       the assumption $q \neq \frac{n+2}{n-2}$ is necessary.
     Consider instead the equation in $B_1\setminus  \{0\}.$ There exists  a nonnegative radial singular solution
      with the blow up rate $|x|^{-\frac{n-2}{2}}$ near the origin.
      Therefore, in Corollary~\ref{t:removal}, the condition $V \in L^{\frac{2n}{n-2}}(B_1)$ is sharp.

 Example 2:  Consider  $\Delta u  + u^q= 0$ in $ \mathbb{R}^n.$   For $q > \frac{n+2}{n-2},$    
   there exists  a   solution asymptotic to $|x|^{-\frac{2}{q -1}}$ at the infinity; see  \cite{KP97}.
  Hence, in Theorem~\ref{t:asymp}, the conditions $q = \frac{2p +n}{n}$ and $V \in L^p$ are sharp.
  Moreover, in Theorem~\ref{t:unique}, the condition   $q = \frac{2p +n}{n}$ is also sharp.

%-------------------------------------------------------------------------------------------------------------------------------------------------------------------
%-------------------------------------------------------------------------------------------------------------------------------------------------------------------
\section{Preliminaries} \label{s:0}
 We collect some standard results in elliptic regularity theory and ordinary differential equations. 
  Lemma~\ref{l:1}-\ref{l:3} follow  by an argument similar to  \cite{BKN89}, section 4.
  
   Let $C_s$ be the Sobolev constant and $\gc = \frac{n}{n-2}.$
  Suppose that the  nonnegative function  $u \in C^{0, 1}$ satisfies $\Delta u + C_0 u^q \geq 0$ weakly in the sense that 
   $$\int (- \langle \gra u, \gra \phi \rangle  + C_0 u^q \phi) dx \geq 0$$
  for all $0 \leq \phi \in C^{\infty}_o.$
  Let $\varphi \geq 0$ be  a function with compact support and $s> 1.$
 Then by Cauchy inequality,
  \begin{eqnarray*}
  \int \varphi^2 u^{q+s-1}dx   &\geq&   C_0^{-1} \int (\frac{4 (s-1)}{s^{2}}   |\varphi \gra u^{\frac{s}{2}}|^2 + 
   \frac{4}{s} \varphi u^{\frac{s}{2}} \langle \gra \varphi,  \gra u^{\frac{s}{2}} \rangle) dx\\
  &\geq&  C_0^{-1} \int (\frac{2}{s^{2}} (s-1)  |\varphi \gra u^{\frac{s}{2}}|^2 - \frac{2}{(s-1)} |\gra \varphi|^2 u^s) dx.
  \end{eqnarray*}
 By Sobolev inequality, we have 
  \beq \label{i:**}
  (\int (\varphi^2 u^s)^{\gc} dx )^{\frac{1}{\gc}} \leq C \int (\frac{s^2 C_0}{2 (s-1)} \varphi^2 u^{q+s-1} + (1+ \frac{s^2}{(s-1)^2}) |\gra \varphi|^2 u^s) dx, 
  \eeq where $C= C(n, C_s, C_0).$ 

 In  Lemma~\ref{l:1}-\ref{l:3}, $u$ is a $C^{0,1}$ function.
   \begin{lemma} \label{l:1} Let  $p> 1$ and $q = \frac{2p+n}{n}.$ 
  Suppose that the nonnegative function   $u \in L^p (B_r)$     
      satisfies $\Delta u + C_0  u^q  \geq 0$ weakly in $B_r.$    
     Then  there exists $\epsilon >0$ such that if $\int_{B_r} u^p dx < \epsilon,$ then
   $\sup_{B_{\frac{r}{2}}} u \leq C r^{-\frac{n}{p}} \|u\|_{L^p (B_r)},$
   where $C=C(n, p, C_s, C_0).$
  \end{lemma} 
    \bpf
   Let $s = p$ in  (\ref{i:**}). Then
   \begin{eqnarray*}
     (\int (\varphi^2 u^p)^{\gc} dx)^{\frac{1}{\gc}} &\leq& C \int (u^{q-1} (\varphi^2 u^p)  +  |\gra \varphi|^2 u^p) dx\\   
        &\leq &  C  (\int _{\{supp\;  \varphi\}} u^p dx)^{\frac{2}{n}} (\int   (\varphi^2 u^p)^{\gc} dx)^{\frac{1}{\gc}} + C \int |\gra \varphi|^2 u^p dx.
     \end{eqnarray*}
     $\varphi$ is chosen to be a cutoff function  such that $\varphi =1$ in $B_{r/2}$ and  $\varphi =0$ outside $B_r$
     with $|\gra \varphi| \leq C r^{-1}.$
     We get
     $$  (\int_{B_{r/2}}  u^{p\gc} dx)^{\frac{1}{\gc}}  \leq \frac{C}{r^2} \int_{B_r} u^p dx.$$
       Choose a sequence $r_k = (2^{-1} + 2^{-k}) r.$
     Apply  (and rescale) the above inequality for $B_{r_k}$ and $B_{r_{k+1}}$ with $p_k = p \gc^{k-1}.$ 
     By Moser iteration, we have
       $\sup_{B_{\frac{r}{2}}} u \leq C r^{-\frac{n}{p}} \|u\|_{L^p (B_r)}.$
  \epf

     \begin{lemma} \label{l:2}  Let  $p> \frac{n}{n-2}$ and $q = \frac{2p+n}{n}.$ 
     Suppose that the nonnegative function   $u \in L^p (\mathbb R^n \setminus B_r)$     
      satisfies $\Delta u + C_0  u^q  \geq 0$ weakly in $\mathbb R^n \setminus B_r.$    
   Then  there exists $\epsilon >0$ such that if $\int_{\mathbb{R}^n \setminus B_r} u^p < \epsilon,$ then $u = O(|x|^{-\gl})$ for all $\gl < n-2$
   as $|x| \rightarrow \infty.$
  \end{lemma}
  
    \bpf
    By Lemma~\ref{l:1}, $u = O(|x|^{-n/p}).$ Let $s = p \frac{n-2}{n} > 1$ in (\ref{i:**}). Then
     $$
     (\int \varphi^{2\gc} u^p dx)^{\frac{1}{\gc}}
          \leq C  (\int _{\{supp \; \varphi\}} u^p dx)^{\frac{2}{n}} (\int  (\varphi^2 u^{p \frac{n-2}{n}})^{\gc} dx)^{\frac{1}{\gc}} +  C \int |\gra \varphi|^2 u^{p \frac{n-2}{n}} dx. 
    $$
   $\varphi$ is chosen to be  a cutoff function such that $\varphi =1$ in $B_{r'} \setminus B_{2r}$ and  $\varphi = 0$ outside
   $B_{2 r'} \setminus B_r$ with $|\gra \varphi| \leq C(1/r + 1/r').$ Let $r' \rightarrow \infty.$ Then
    $$  (\int \varphi^{2\gc} u^p dx)^{\frac{1}{\gc}} \leq C (\int |\gra \varphi|^n dx)^{\frac{2}{n}} (\int_{\{supp\;  \gra \varphi\}} u^p dx)^{\frac{1}{\gc}}.$$
  And thus, 
    $$ (\int_{\mathbb R^n \setminus B_{2r}}  u^p dx)^{\frac{1}{\gc}}  \leq C (\int_{B_{2r} \setminus B_{r}} u^p dx)^{\frac{1}{\gc}}.$$  
  This gives $\int_{\mathbb R^n \setminus B_{r}}  u^p = O(r^{- \gd})$ for some small $\gd > 0.$
  Therefore, by Lemma~\ref{l:1}, $u = O(|x|^{-\frac{n}{p}- \frac{\gd}{p}}).$ 
  Let $\gl_0 = \sup \{\gl : u = O(|x|^{-\gl}) \}.$ By iteration and  a contradiction argument, we get that $\gl_0 = n-2.$ 
  \epf 

   Suppose that $h \geq 0$ is a $C^0$ function.
   The  nonnegative function  $u \in C^{0, 1}$ satisfies $\Delta u + C_0 h u \geq 0$ weakly if 
   $$\int (- \langle \gra u, \gra \phi \rangle  + C_0 h u \phi) dx \geq 0$$
  for all $0 \leq \phi \in C^{\infty}_o.$
     \begin{lemma} \label{l:3} Let  $p> 1$ and  $t > \frac{n}{2}.$ Suppose that  the nonnegative function  $h \in L^t (B_r)$ satisfies
      $\int_{B_r} h^t dx \leq \frac{C_1}{r^{2t-n}}.$  Suppose also that the nonnegative function   $u \in L^p (B_r)$     
      satisfies $\Delta u + C_0 h u \geq 0$ weakly in $B_r.$ 
     Then $\sup_{B_{\frac{r}{2}}} u \leq C r^{-\frac{n}{p}} \|u\|_{L^p (B_r)},$
     where $C= C(n, p, C_s, C_0, C_1).$
    \end{lemma}
    \bpf
   The proof is by standard Moser iteration. See  Morrey \cite{Morrey}. 
  \epf

The following is a basic result in ordinary differential equations; see  \cite{Chen08}.
\begin{lemma} \label{l:ode}
  Suppose that $f(r) \geq 0$ satisfies $f(r) \leq - \frac{r}{a} f'(r) + C_2 r^{-b}$ for some $a, b > 0.$
  
  (a) $a \neq b.$ Then there exists a constant $C_3$ such that
      $$f(r) \leq C_3 r^{-a} + \frac{a \,C_2}{a-b} r^{-b}.$$ Therefore, $f(r) = O(r^{-\min \{a, b\} })$ as $r \rightarrow \infty.$ 
  
  (b) $a = b.$ Then there exists a constant $C_3$ such that
       $$f(r) \leq C_3 r^{-a} + a \,C_2 r^{-a}\ln r.$$ Therefore, $f(r) = O(r^{-a} \ln r)$ as $r \rightarrow \infty.$ 
 \end{lemma}

%------------------------------------------------------------------------------------------------------------------------------------------------------------------- %-------------------------------------------------------------------------------------------------------------------------------------------------------------------

\section{Proof of Theorem~\ref{t:asymp}} \label{s:1}

   \bpf [Proof of Theorem~\ref{t:asymp}] 
    We first derive a version of Pohozaev's identity for non-variational systems. 
Let $\Omega$ be a domain in $\mathbb R^n$ and $N$ be the unit outer normal on $\de \Omega.$  We will perform integration by parts repeatedly. 
  \begin{align}
   & \int_{\Omega} \sum_{k,l} f^k(V) x_l D_l V_k dx = \int_{\Omega} \sum_{j, k,l} A_{kj} \Delta V_j x_l D_l V_k dx\notag\\
    & 
    =   \int_{\Omega} - \sum_{i,j, k, l} A_{kj}  D_iV_j D_i (x_l D_lV_k) dx + \int_{\de \Omega} \sum_{i,j, k, l} A_{kj} D_i V_j x_l D_lV_k N_i  dS\notag\\
     &=   \int_{\Omega} (- \sum_{i,j, k} A_{kj} D_i V_j D_iV_k - \sum_{i,j, k, l} D_l( A_{kj} D_i V_j D_i V_k) \frac{x_l}{2}) dx +
     \int_{\de \Omega} \sum_{i,j, k, l} A_{kj} D_i V_j x_l D_lV_k N_i dS\notag\\
    &=   (\frac{n}{2} - 1) \int_{\Omega}   \sum_{i, j, k}    A_{kj} D_i V_j D_iV_k dx - \int_{\de \Omega} \frac{1}{2}\sum_{i, j,k, l}  A_{kj} D_i V_j D_i V_k x_l N_l dS\notag\\
     &  + \int_{\de \Omega} \sum_{i,j, k, l} A_{kj} D_i V_j x_l D_lV_k N_i dS. \label{e:pohoz}
  \end{align}
   It is worth mentioning that $x_l D_l$ is a conformal killing field in $\mathbb{R}^n.$
  
   As we explained in the introduction, without loss of generality we may assume the equation is of the diagonal form, i.e., 
    \beq \label{e:V}
   \Delta V_i = f^i (V).
   \eeq 
   Note that $|V|$ and $|\gra V|$ are $C^{0,1}$ functions.
  By (\ref{e:V}) and (A1)- (A2),
  we have
  $$\begin{array}{l}
   \Delta |V| \geq - C |V|^q; \\
    \Delta |\gra V| \geq - C |V|^{q-1} |\gra V|
 \end{array}
  $$ weakly. 
  % If $\tilde p < p,$ then $V \in L^{\tilde p} (\mathbb{R}^n - K)$ implies
  %  $V \in L^p(\mathbb{R}^n - K).$ Without loss of generality, we may assume $q= \frac{2p +n}{n}.$
 Since $V \in L^p (\mathbb R^n \setminus K),$ there exists  a large number $R$ such that $\int_{\mathbb R^n \setminus B_R} |V|^p dx < \epsilon,$
 where $\epsilon$ is as in Lemma~\ref{l:1}. 
 Applying Lemma~\ref{l:1} to  $B_r(x_0)$  where $|x_0| \geq 2r \geq 2 R,$ we get $|V|= O(|x|^{- \frac{n}{p}}).$

 Case 1. If $\frac{n+2}{n} < q \leq \frac{n}{n-2}$ (or equivalently, $1< p \leq \frac{n}{n-2}$), then $\frac{n}{p} \geq n-2.$ By Lemma~\ref{l:1}, we have $|V|= O(|x|^{-n/p}).$ Let $\varphi$ be a cutoff function such that $\varphi=1$ in $B_r$ and $\varphi=0$ outside $B_{2r}$ with
   $|\gra \varphi| \leq Cr^{-1}.$  
  Applying $\varphi V_i$ to (\ref{e:V}) and integrating gives
              $$\int_{B_r(x_0)} |\gra V|^2 dx \leq C \int_{B_{2r}(x_0)} |V|^{q+1} dx + \frac{C}{r^2} \int_{B_{2r} (x_0)} |V|^2 dx = O(r^{n-2 - \frac{2n}{p}}) \leq O(r^{-n +2}),$$ 
                where $|x_0| \geq 2r \gg 1.$    
                        By Lemma~\ref{l:3} with $h = |V|^{q-1},$ we obtain $|\gra V|= O(|x|^{-(n-1)})$ and thus $|V| = O(|x|^{-(n-2)}).$

 Case 2. If $\frac{n}{n-2} < q$ (or equivalently $p > \frac{n}{n-2}$), by Lemma~\ref{l:2}, $|V|= O(|x|^{-\gl})$ for all $\gl < n-2.$
                Therefore, $$\int_{B_r(x_0)} |\gra V|^2 dx \leq C \int_{B_{2r}(x_0)} |V|^{q+1} dx + \frac{C}{r^2} \int_{B_{2r} (x_0)} |V|^2 dx = O(r^{n-2 - 2 \gl}),$$ 
                where $|x_0| \geq 2r \gg 1.$  Moreover, 
                $|V| \in L^{p'}$ for all $p' > \frac{n}{n-2}.$ Choose $p' < p$  close to $\frac{n}{n-2}.$ 
                Hence, $q > \frac{2p' +n}{n}.$ We can then find $q' > \frac{n}{2}$ such that
                $$\int_{B_r(x_0)} (|V|^{q-1})^{q'} dx \leq \frac{C}{r^{2q'-n}}$$  where $|x_0| \geq 2r \gg 1.$This is possible because $\gl$ is close to $n-2.$
                By Lemma~\ref{l:3}, we obtain
           $$
                 \sup_{B_{\frac{r}{2} (x_0)}} |\gra V| \leq \frac{C}{r^{\frac{n}{2}}} \|\gra V\|_{L^2(B_r (x_0))}= O(r^{- \gl -1}),               
            $$    where $|x_0| \geq 2r \gg 1.$
   
      Let $\Omega= B_R \setminus B_r$ in (\ref{e:pohoz}). We have
    \begin{align}
   \int_{\Omega} \sum_{k,l} f^k(V) x_l D_l V_k dx &=   (\frac{n}{2} - 1) \int_{\Omega}       |\gra V|^2 dx 
        - \int_{\de \Omega} \frac{1}{2}  \sum_l |\gra V|^2 x_l N_l dS\notag\\
     & + \int_{\de \Omega} \sum_{i,j, l}  D_i V_j x_l D_lV_j N_i dS. \label{i:case2}
    \end{align}
   Note that  
     $$\lim_{R \rightarrow \infty} \int_{S_R} R |\gra V|^2 dS = \lim_{R \rightarrow \infty} O(R^{-2 \gl - 2+n})
   = 0.$$ 
    Let  $R \rightarrow \infty$ in (\ref{i:case2}). Then there is no boundary term coming from the infinity. 
    We can choose
    $\Omega = \mathbb{R}^n \setminus B_r.$  The boundary terms only occur on $S_r.$ On $\de \Omega,$ $N = - \frac{x}{r}.$ 
   Hence,
    \begin{eqnarray*}
   \int_{ \mathbb{R}^n \setminus B_r} \sum_{k,l} f^k(V) x_l D_l V_k dx &=&   (\frac{n}{2} - 1) \int_{ \mathbb{R}^n \setminus B_r}     |\gra V|^2 dx 
        + \int_{S_r} \frac{r}{2}   |\gra V|^2 dS
      - r \int_{S_r}  |\gra_N V|^2 dS\\
      &\geq&  (\frac{n}{2} - 1) \int_{ \mathbb{R}^n \setminus B_r}     |\gra V|^2 dx 
        - \int_{S_r} \frac{r}{2}   |\gra V|^2 dS.
   \end{eqnarray*}
    Let $G(r) := \int_{ \mathbb{R}^n \setminus B_r}   |\gra V|^2 dx.$ Since $G'(r) = - \int_{S_r}   |\gra V|^2 dS,$
     the previous formula becomes
     $$ G(r)
      \leq   -\frac{r}{n-2} G'(r)  +  \frac{2}{n-2} \int_{ \mathbb{R}^n \setminus B_r} \sum_{k,l} f^k(V) x_l D_l V_k dx.$$
    
     The key idea is to compare the size of $G(r)$ to that of $G'(r)$.    The coefficient in front of $G'(r)$  plays an important role. 
     Here is the only place we use the condition of $|f^i|.$ 
     We have
      $$ \int_{ \mathbb{R}^n \setminus B_r} \sum_{k,l} f^k(V) x_l D_l V_k dx \leq \int_{ \mathbb{R}^n \setminus B_r}|V|^q |x| |\gra V| dx = O(r^{- \gl(q+1)+ n}).$$
    Thus,      
    $$G(r) \leq - \frac{r}{n-2} G'(r) +C r^{- \gl(q+1)+ n}.$$
    Since $q > \frac{n}{n-2}$ and $\gl$ is close to $n-2,$ we have $\gl(q+1)- n> n-2.$
   By  Lemma~\ref{l:ode},  this implies $G(r) = O(r^{-(n-2)}).$ 
  By Sobolev inequality,  we get 
   $$\int_{B_{2r} \setminus B_r} |V|^{\frac{2n}{n-2}} dx = O(r^{- n}).$$ 
   Finally, by Lemma~\ref{l:1} and \ref{l:3} we obtain
  $|V| = O(|x|^{-(n-2)})$ and $|\gra V| = O(|x|^{- (n-1)}).$ 
  \epf

 %-------------------------------------------------------------------------------------------------------------------------------------------------------------------
\section{Proofs of Corollary~\ref{t:removal} and Theorem~\ref{t:unique}} \label{s:2}

  \bpf [Proof of Corollary~\ref{t:removal}]  
    Since the equation is invariant under inversion, we  transform the solution to $\mathbb{R}^n \setminus B_1$ and apply Theorem~\ref{t:asymp}. 
   
     Let $y = \frac{x}{|x|^2}.$ Define $U_i (y)= \frac{1}{|y|^{n-2}} V_i (\frac{y}{|y|^2}).$
     This is called the Kelvin transform with the property that
    $$\Delta_y U_i(y) = \frac{1}{|y|^{n+2}} \Delta_x V_i (x).$$     
    This can also be viewed as the  conformal change formula of the conformal Laplacian with zero scalar curvature. 
     Therefore, $U_i (y)$ satisfies
      $$\sum_j A_{ij} \Delta_y U_i (y) = \frac{1}{|y|^{n+2}} f^i (|y|^{n-2} U(y)) = f^i (U(y))$$
    in $\mathbb{R}^n \setminus B_1,$ where we use that $f^i$ is homogeneous of degree $\frac{n+2}{n-2}.$ Moreover, 
    \begin{eqnarray*}
    \int_{\mathbb{R}^n \setminus B_1} |U|^{\frac{2n}{n-2}} dy &=& \int_{\mathbb{R}^n \setminus B_1} (|V| |y|^{-n+2})^{\frac{2n}{n-2}} dy = \int_{B_1\setminus \{0\}}  (|V| |x|^{n-2})^{\frac{2n}{n-2}} |x|^{-2n} dx\\
    &=&  \int_{B_1\setminus \{0\}} |V|^{\frac{2n}{n-2}} dx  < + \infty.
     \end{eqnarray*}    
     Now we apply Theorem~\ref{t:asymp} with $p = \frac{2n}{n-2}$ and $q = \frac{n+2}{n-2}.$ We get $|U| = O(|y|^{- (n-2)})$ and $|\gra U| = O(|y|^{-(n-1)}).$
     Hence,  $|V| = O(1)$ and $|\gra V| = O(|x|^{-1}).$ As a result, $V \in L^{\infty}(B_1)$ and $\gra V \in L^p(B_1)$ for all $p< n.$
     
     We show that $V$ is a weak solution to (\ref{e:system})  in $B_1.$ Let $\varphi \in H^1_0 (B_1, \mathbb{R}^m).$ Let  $\eta_k(|x|)$ be a compactly supported function in $B_1 \setminus \{0\}$
     such that $\eta_k \rightarrow 1$ a.e. in $B_1$ and $\|\eta_k \|_{L^n(B_1)} \rightarrow 0$ as $k \rightarrow \infty.$ (Such functions were used by Serrin \cite{Serrin64}.)
     Then 
     $$
     \int_{B_1} \eta_k \sum_{i,j, l} A_{ij} D_l \varphi_j D_l V_i dx= \int_{B_1} - \sum_i f^i(V)  \varphi_i \eta_k dx- \int_{B_1}  \sum_{i,j, l}   D_l \eta_k  A_{ij} \varphi_j D_l V_i dx. 
    $$
    The last term can be estimated as follows.
    $$\vert \int_{B_1}  \sum_{i,j, l}   D_l \eta_k  A_{ij} \varphi_j D_l V_i  dx \vert \leq C  \|\varphi\|_{L^{\frac{2n}{n-2}}(B_1)} \|\gra V\|_{L^2(B_1)} \|\eta_k\|_{L^n(B_1)}
     \leq C \|\eta_k\|_{L^n(B_1)} \rightarrow 0$$
    as $k \rightarrow \infty.$ Hence, in the limit
     $$ \int_{B_1} \sum_{i,j, l} A_{ij} D_l \varphi_j D_l V_i dx = \int_{B_1} - \sum_i f^i(V)  \varphi_i dx.$$
    Thus, $V$ is a weak solution in $B_1$. It follows by elliptic regularity that $V \in C^{\infty} (B_1).$
  \epf

  \bpf [Proof of Theorem~\ref{t:unique}] 
 Let $\Omega = B_R.$ Therefore, $N = \frac{x}{R}$  in (\ref{e:pohoz}). We get
   \begin{eqnarray*}
   \int_{B_R} \sum_{k,l} f^k(V) x_l D_l V_k dx &=&   (\frac{n}{2} - 1) \int_{B_R}   \sum_{i, j, k}    A_{kj} D_i V_j D_iV_k dx\\
   &-& \int_{S_R} \frac{R}{2}\sum_{i, j,k}  A_{kj} D_i V_j D_i V_k dS
       + R \int_{S_R} \sum_{j, k} A_{kj} D_N V_j D_NV_k dS.
    \end{eqnarray*}  Since $f^k = \frac{\de P}{\de V_k},$ we have
     \begin{align}
   \int _{B_R} -n P(V) dx
   &=   (\frac{n}{2} - 1) \int_{B_R}   \sum_{i, j, k}    A_{kj} D_i V_j D_iV_k dx- \int_{S_R} \frac{R}{2}\sum_{i, j,k}  A_{kj} D_i V_j D_i V_k dS \notag\\
      &  + R \int_{S_R} \sum_{j, k} A_{kj} D_N V_j D_NV_k dS- \int_{S_R} R P(V) dS. \label{e:1}
    \end{align}    
    On the other hand, we also have
     \begin{align}
   \int _{B_R} (q+1) P(V) dx
   &= \int_{B_R} \sum_k \frac{\de P}{\de V_k} V_k dx \notag\\
   &=  - \int_{B_R}   \sum_{i, j, k}    A_{kj} D_i V_j D_iV_k dx
       +  \int_{S_R} \sum_{j, k} A_{kj} D_N V_j V_k dS,  \label{e:2}
    \end{align}    where we use the Euler formula for homogeneous functions.

  Case 1.    $n \geq 4.$
  
   By Theorem~\ref{t:asymp},  when $R \rightarrow \infty,$ (\ref{e:1}) becomes
    \begin{eqnarray*}
      \int _{B_R} -n P(V) dx
  & =&   (\frac{n}{2} - 1) \int_{B_R}   \sum_{i, j, k}    A_{kj} D_i V_j D_iV_k dx + O(R^{-(n-2)})+ O(R^{- (q+1)(n-2) + n})\\
  &=& (\frac{n}{2} - 1) \int_{B_R}   \sum_{i, j, k}    A_{kj} D_i V_j D_iV_k dx + o(1),
   \end{eqnarray*} where we use conditions on  $p, q$ and $n \geq 4$ to get $(q+1)(n-2) - n > 0.$
   Similarly, (\ref{e:2}) gives 
  $$  \int _{B_R} (q+1) P(V) dx
  =  - \int_{B_R}   \sum_{i, j, k}    A_{kj} D_i V_j D_iV_k dx+ O(R^{- (n-2)}).
      $$
     Combining these two formulas and noting that $q+1 \neq \frac{2n}{n-2},$ we finally arrive at
     $$ \int_{B_R}   \sum_{i, j, k}    A_{kj} D_i V_j D_iV_k dx = o(1).$$
      Since $A$ is positive definite, we have $|\gra V| \equiv 0$ and hence $V \equiv 0.$ 
      
    Case 2. $n=3.$  
    
    Note that  $\sup |V| \leq \frac{C}{|x|^{n/p}} \|V\|_{L^p}.$  Combining this fact with Theorem~\ref{t:asymp}, we have $|V| = O(|x|^{- \gl}),$
    where $\gl = \max \{1, \frac{3}{p}\}.$ Therefore, 
    $$\gl (q+1)-3 \geq \max \{q-2, \frac{3}{p} (q+1) - 3\} \geq \max \{-1+ \frac{2p}{3}, -1 + \frac{6}{p}\} > 0.$$
    Then (\ref{e:1}) becomes
      \begin{eqnarray*}
      \int _{B_R} -3 P(V) dx
  & =&   (\frac{3}{2} - 1) \int_{B_R}   \sum_{i, j, k}    A_{kj} D_i V_j D_iV_k dx + O(R^{-1})+ O(R^{- \gl (q+1) + 3})\\
  &=& (\frac{3}{2} - 1) \int_{B_R}   \sum_{i, j, k}    A_{kj} D_i V_j D_iV_k dx + o(1)
   \end{eqnarray*}  as in Case 1.
    The rest of proof is the same as in Case 1.
  \epf

 \textbf{Acknowledgments:}  The author would like to thank Craig Evans for drawing her attention to the reaction-diffusion
   type elliptic systems. She is  grateful for many of his suggestions. The author appreciates James Serrin for his comments and especially
     suggestions for the references, which help improve the presentation of the work.   The author wishes to thank Haim Brezis for comments 
     and interests. Finally, the author would like to thank the referee for careful reading and for valuable suggestions. 

  While the author was preparing for the manuscript, she was supported by NSF grant DMS-0635607.

 % Department of Mathematics, University of California, Berkeley, CA
 %\par
% Email address: \textsf{sophie@math.berkeley.edu}

 Institute for Advanced Study, Princeton, NJ
 \par
 Email address: \textsf{sophie@math.ias.edu}

\end{document}